\input amstex
\input amsppt.sty
\magnification1200
\vsize=23.5truecm
\hsize=16.5truecm
\vcorrection{-10truemm}
\NoBlackBoxes

\def\Ama{A_{\max}}

\def\wA{\widetilde A}

\comment

\def\comega{\overline\Omega }
\def\simto{\overset\sim\to\rightarrow}
\def\ang#1{\langle {#1} \rangle}

\def\rp{ \Bbb R_+}

\define\tr{\operatorname{tr}}
\endcomment

\input gerd3.sty

\document
\topmatter
\title
Spectral asymptotics for Robin problems with a discontinuous coefficient
\endtitle
\author Gerd Grubb \endauthor
\affil
{Department of Mathematical Sciences, Copenhagen University,
Universitetsparken 5, DK-2100 Copenhagen, Denmark.
E-mail {\tt grubb\@math.ku.dk}}\endaffil

\abstract
The spectral behavior of the difference between the resolvents of two
realizations $\wA_1$ and $\wA_2$ of a second-order strongly 
elliptic symmetric differential operator $A$, defined by 
different Robin conditions $\nu u=b_1\gamma_0u$ and $\nu
u=b_2\gamma_0u$, can in the case where all coefficients are $C^\infty$ be
determined by use of a general result by the author in 1984 on singular Green
operators. We here treat the problem for nonsmooth $b_i$. Using a
Krein resolvent formula, we show that if
$b_1$ and $b_2$ are in $L_\infty$, the s-numbers $s_j$ of $(\wA_1   -\lambda
)^{-1}-(\wA_2   -\lambda )^{-1}$  satisfy $s_j j^{3/(n-1)}\le
C$ for all $j$; this improves a recent result for $A=-\Delta $ 
by Behrndt et al., that
$\sum_js_j ^p<\infty$ for $p>(n-1)/3$. 
A sharper estimate is obtained when $b_1$ and $b_2$ are in
$C^\varepsilon$ for some $\varepsilon >0$, with jumps at a smooth
hypersurface, namely that $s_j j^{3/(n-1)}\to c$ for $j\to \infty$, with a
constant $c$ defined from the principal symbol of $A$ and $b_2-b_1$. 

As an auxiliary result we show that the usual principal spectral asymptotic estimate for
pseudodifferential operators of negative order on a closed manifold 
extends to products of pseudodifferential operators 
interspersed with piecewise continuous functions.
\endabstract
\dedicatory Dedicated to the memory of M.\ Sh.\ Birman (1928--2009) \enddedicatory
\endtopmatter
\rightheadtext {Robin problems}

\subhead  Introduction \endsubhead

Consider a second-order strongly elliptic symmetric operator $$
A=-\sum_{j,k=1}^n\partial_j(a_{jk}\partial_ku)+a_0u\tag0.1
$$ on a
bounded smooth domain $\Omega \subset{\Bbb R}^n$, and denote by
$A_\gamma $, $A_{  \nu } $, resp.\ $\wA$, the realizations in $L_2(\Omega
)$ defined by the Dirichlet condition $\gamma _0u=0$, the Neumann
condition
$\nu u=0$, resp.\ a Robin condition $\nu u-b\gamma _0u=0$ with $b$ real. Here $\gamma
_0u=u|_{\partial\Omega }$, and 
$$
\nu u = \sum_{j,k=1}^n n_j \gamma_0 (a_{jk}\partial_{ k} u),\tag0.2
$$
the conormal derivative, 
with $\vec n=(n_1,\dots, n_n)$ denoting the
interior normal to $\partial\Omega $.
It is a classical result of Birman \cite{B62}, shown also for exterior
domains, that the difference between the resolvents of the Robin
realization and the Dirichlet realization is compact and has the
spectral behavior, for large negative
$\lambda$,
$$
s_j\big((\wA-\lambda )^{-1}-(A_\gamma    -\lambda )^{-1}\big)j^{2/(n-1)}\le C
\text{ for all }j ;\tag0.3
$$
here $s_j(T)$ denotes the $j$-th eigenvalue of $(T^*T)^{\frac12}$ (the
$j$-th s-number or singular value of $T$), counted
with multiplicities. This was shown assuming merely that $b\in L_\infty
(\partial\Omega) $. 

For the situation where all coefficients are
$C^\infty $, the estimate was later improved to an asymptotic estimate
$$
s_j\big((\wA-\lambda )^{-1}-(A_\gamma    -\lambda )^{-1}\big)j^{2/(n-1)}\to c
\text{ for }j\to\infty  ;\tag0.4
$$
this follows from Grubb \cite{G74}, Sect.\ 8 (with generalizations to higher-order operators),
and Birman and Solomiak \cite{BS80} (including exterior domains).

The paper \cite{G84} gave tools to extend (0.4) to
nonselfadjoint situations (also for exterior domains by a cutoff technique), by showing that
for any singular Green operator $G$ on $\Omega $ of order $-t<0$ and class 0,
$$
s_j(G)j^{t/(n-1)}\to c(g^0) \text{ for }j\to \infty ;\tag0.5
$$
here $G$ belongs to the calculus of pseudodifferential boundary operators,   
introduced by Boutet de Monvel \cite{B71} and further developed in
\cite{G84}, \cite{G96}; $c(g^0)$ is a constant derived from the
principal symbol $g^0$. In fact, the resolvent difference 
in (0.4) is a singular Green operator of order $-2$ and class 0, when all
coefficients are $C^\infty $.

Considering another resolvent difference, J. Behrndt, M. Langer, 
I. Lobanov,
V. Lotoreichik and I. Popov showed in a recent paper \cite{BLLLP10}, on the
basis of a theory of quasi-boundary triples by J. Behrndt
and M. Langer \cite{BL07}, that when $A=-\Delta $ (hence $\nu
u=\gamma _1u=\sum_jn_j\gamma_0\partial_{ j}u$) and
$b$ is a real function in $L_\infty
(\partial\Omega )$, the
difference between the resolvent of $\wA$ and the resolvent of the
Neumann realization $A_{\nu  } $ satisfies an estimate with 2 replaced
by 3, for $\lambda $ in the
intersection of resolvent sets $\varrho (\wA )\cap\varrho (A_\nu  )$:
$$
(\wA-\lambda )^{-1}-(A_\nu -\lambda )^{-1}\in \Cal C_p \text{ for }p>3/(n-1);\tag0.6
$$
here $\Cal C_p$ denotes  the space of compact operators $T$ with singular value sequences $(s_j(T))_{j\in{\Bbb N}}\in
\ell_p$; the Schatten class of order $p$. 
(Besides real $b$, also functions with a fixed sign on
$\operatorname{Im}b$ were treated.) 

In the case where $b\in C^\infty
(\partial\Omega )$, this follows from (0.5) since $(\wA-\lambda )^{-1}-(A_{\nu  }
-\lambda )^{-1}$ is a singular Green operator of order $-3$, 
leading to the stronger estimate:
$$
s_j\big((\wA-\lambda )^{-1}-(A_{\nu  }   -\lambda )^{-1}\big)j^{3/(n-1)}\to c
\text{ for }j\to \infty ;\tag0.7
$$
this was noted also in
\cite{G10a}, Cor.\ 8.4 and Ex.\ 8.5.

The result of \cite{BLLLP10} is more general by treating nonsmooth
$b$. Their main result Theorem 3.5 is proved in a formulation where
the boundary condition is $\gamma _0u=\Theta \gamma _1u$ for an
operator $\Theta $, but it is explained in their
Remark 3.7 how this can be made to include general conditions $\gamma
_1u=b\gamma _0u$ by use of the relations point of view of \cite{BL07};
more on this in \cite{BL10}.

The main purpose of the present paper is to show spectral asymptotics
estimates as in (0.7) for  nonsmooth $b$. First we show (in Theorem 2.2) that an upper bound
$$
s_j\big((\wA-\lambda )^{-1}-(A_\nu   -\lambda )^{-1}\big)j^{3/(n-1)}\le C
\text{ for all }j ,\tag0.8
$$
holds for any complex $b\in L_\infty
(\partial\Omega )$; this implies (0.6), for general $A$ as in (0.1). As a
corollary, a similar estimate holds for  $(\wA_1   -\lambda
)^{-1}-(\wA_2   -\lambda )^{-1}$, when the $\wA_i$ are defined by
boundary conditions $\nu u=b_i\gamma _0u$ with $b_i\in L_\infty
(\partial\Omega )$.
 
Next, we show (in Theorems 2.5 and 3.4) that asymptotic estimates hold when the $b_i$ are
piecewise slightly better that continuous. Since asymptotic estimates are not additive, we must
aim directly for $(\wA_1   -\lambda
)^{-1}-(\wA_2   -\lambda )^{-1}$. 

\proclaim{Theorem 0.1} Assume that $b_1$, $b_2$ and $b_2-b_1$ are
piecewise $C^\varepsilon $ on $\partial\Omega $ for some $\varepsilon >0$, having jumps at
$C^\infty $ hypersurfaces. Then
$$
s_j\big((\wA_1   -\lambda
)^{-1}-(\wA_2   -\lambda )^{-1})\big)j^{3/(n-1)}\to c
\text{ for }j\to \infty ,\tag0.9
$$
where $c$ is defined from  the principal symbol of $A$ and $b_2-b_1$.
\endproclaim

It suffices in fact that  $b_1$, $b_2$ and $b_2-b_1$ are piecewise in $H^r_p(\partial\Omega
)$ for some $r>0$ and $p>(n-1)/r$, see the details below.

For the proof of (0.8) the method is, as in \cite{BLLLP10}, an
application of functional analysis, building on a general
theory of extensions (here Grubb \cite{G68}) 
together with known facts on
elliptic boundary value problems. The proof of (0.9) in
the nonsmooth situations draws on methods and results for 
pseudodifferential 
boundary operators in \cite{G84} and a result on 
restricted kernels of pseudodifferential operators by Laptev \cite{L77, L81}. 

As an auxiliary result
of independent interest we show (Theorem 3.3) that a product of classical
pseudodifferential operators of negative order on a closed manifold,
interspersed with piecewise continuous functions having jumps at a
smooth hypersurface, has a principal spectral asymptotics estimate as
in the smooth case. Moreover, we extend (0.4) to $b\in L_\infty
(\partial\Omega )$ (Theorem 2.4).

Some spectral
estimates for resolvent differences in interior and exterior domains
have been described recently by Malamud in [M10], and spectral
asymptotics have been shown in [G10]; both papers treat higher-order operators
but do not aim for the special bounds obtained here.
Let us also mention that we do not here address the question of 
nonsmooth domains, as e.g.\ in
Gesztesy and Mitrea \cite{GM09, GM09a, GM10} and Abels, Grubb and Wood
\cite{AGW10}, \cite{G08}, and their references.

\medskip

To keep the paper short, some introductory material
found in other sources will not be repeated here.

The main details of the extension theory \cite{G68}--\cite{G74} have 
been recalled and
explained in several recent papers  \cite{BGW09}, \cite{G08},
\cite{G10a}; resulting Krein-type resolvent formulas are shown in \cite{BGW09}.

Sobolev spaces are recalled in numerous places. The basic facts we
shall need on these and
other function spaces such as Besov and Bessel-potential
spaces, are recalled e.g.\ in \cite{AGW10}, Sect.\ 2.

The calculus of pseudodifferential boundary operators is explained in Boutet de Monvel \cite{B71} and
in \cite{G84}, \cite{G96}, \cite{G09}.

\subhead 1. The Robin realization \endsubhead

Let $\Omega $ be a bounded smooth subset of $\rn$ with boundary $\partial\Omega =\Sigma $, and let 
$$
a(u,v)=\sum _{j,k=1}^n (a_{jk}\partial _{ k}u, \partial _{ j}v)
+(a_0u, v),\tag1.1
$$
be a sesquilinear form with coefficients in $C^\infty (\comega )$ such that
the associated 
second-order operator (0.1)
is formally selfadjoint and strongly elliptic. We assume moreover 
that $a(u,u)$ is real for $u\in H^1(\Omega )$ and (with $c>0$, $k\ge 0$)
$$
a(u,u)\ge c\|u\|^2_1-k\|u\|^2_0,\text{ for }u\in H^1(\Omega ).\tag1.2
$$
This holds if the matrix $(a_{jk}(x))_{j,k=1}^n$ is real, symmetric and
positive definite 
and $a_0(x)$ is real, at each
$x\in\comega$.

Let $b\in L_\infty (\Sigma  )$, and define the sesquilinear form $a_b$ by
$$
a_b(u,v)=a(u,v)+(b\gamma _0u,\gamma _0v)_{L_2(\Sigma )}.\tag1.3
$$
Since $\|\gamma _0u\|^2_{L_2(\Sigma )}\le c'\|u\|^2_{\frac34}\le \varepsilon
\|u\|^2_{1}+C(\varepsilon )\|u\|^2_0$ for any $\varepsilon
$, we infer from (1.2) that
$$
\operatorname{Re}a_b(u,u)\ge c_1\|u\|^2_1-k_1\|u\|^2_0,\text{ for }u\in H^1(\Omega ),\tag1.4
$$
where $c_1<c$ is close to $c$ and $k_1\ge k$ is a large constant.

The sesquilinear form $a_b$ on $V=H^1(\Omega )$ in $H=L_2(\Omega )$
defines a realization $\wA$ of $A$ by Lions' version of the Lax-Milgram
lemma (as recalled e.g.\ in \cite{G09}, Ch.\ 12), with domain
$$
D(\wA)=\{u\in H^1(\Omega )\cap D(\Ama)\mid (Au,v)=a_b(u,v)\text{ for all
}v\in H^1(\Omega )\}.\tag1.5
$$
The operator
$\wA$ is
closed, densely defined with spectrum in a sectorial region in
\linebreak$\{\operatorname{Re}\lambda \ge -k_1\}$, and its adjoint $\wA^*$ is the
analogous operator defined from $$
a_b^*(u,v)=\overline{a(v,u)}+(\overline b \gamma _0u,\gamma
_0v)_{L_2(\Sigma )}.\tag1.6$$
In particular, when $b$ is real, $\wA$ is
selfadjoint. 

It will be useful to observe:

\proclaim{Lemma 1.1} For any small $\theta >0$ there is an
$\alpha \ge 0$ such that the spectrum of $\wA$ is contained in the region
$$
M_{\theta ,\alpha ,k_1}=\{z\in {\Bbb C}\mid |\operatorname{Im}z|\le
\theta (\operatorname{Re}z+\alpha),\; \operatorname{Re}z\ge -k_1 \}.\tag1.7
$$
\endproclaim

\demo{Proof} Let $K=\|\operatorname{Im}b\|_{L_\infty (\Sigma )}$. From
the inequalities for $a_b(u,u)$ we see that for $u\in H^1(\Omega )$,
$$
\aligned
|\operatorname{Im}a_b(u,u)|&=|\operatorname{Im}(b\gamma _0u,\gamma
_0u)|\le K(\varepsilon \|u\|_1^2+C(\varepsilon )\|u\|_0^2)\\
&\le K\varepsilon
c_1^{-1}(\operatorname{Re}a_b(u,u)+k_1\|u\|_0^2)+KC(\varepsilon
)\|u\|_0^2\\
&=K\varepsilon c_1^{-1}\operatorname{Re}a_b(u,u)+
(K\varepsilon c_1^{-1}k_1+KC(\varepsilon ))\|u\|_0^2.
\endaligned
$$
This (together with (1.4)) shows that for $u\ne 0$, $a_b(u,u)/\|u\|_0^2$
has its
values in $M_{\theta ,\alpha ,k_1}$, where $\theta =K\varepsilon
c_1^{-1}$ can be taken  arbitrarily small, $\alpha =K\varepsilon c_1^{-1}k_1+KC(\varepsilon )$. The numerical ranges of $\wA$ and
$\wA^*$ are contained in this set, which then also contains the spectra.
(More details for this kind of argument can be found in \cite{G09},
Sect.\ 12.4.)\qed \enddemo

The Neumann-type boundary operator (0.2)  
enters in the ``halfways Green's formula''
$$
(Au,v)-a(u,v)=(\nu u,\gamma _0v)_{L_2(\Sigma )} ,\tag1.8
$$
for smooth $u$ and $v$. It is known e.g.\ from \cite{LM68} that
$\gamma_1$ and $\nu $ extend to continuous mappings from  $ H^1(\Omega )\cap D(\Ama)$ to
$H^{-\frac12}(\Sigma )$, such that for $u\in  H^1(\Omega )\cap
D(\Ama)$, $v\in H^1(\Omega )$,  (1.8) holds with the scalar product
over $\Sigma $ replaced by the sesquilinear duality between
$H^{-\frac12}(\Sigma )$ and $H^\frac12(\Sigma )$.
Then
$$
(Au,v)-a_b(u,v)=(\nu u,\gamma _0v)_{H^{-\frac12}(\Sigma ), H^\frac12(\Sigma )}-(b\gamma _0u,\gamma _0u)_{L_2(\Sigma )},\tag1.9
$$
and hence
$$
D(\wA)=
\{u\in H^1(\Omega )\cap D(\Ama)\mid \nu  u=b\gamma _0u\text{ in }H^{-\frac12}(\Sigma )\},
$$
representing the Robin condition $\nu u-b\gamma _0u=0$.

For $b=0$, the condition is $\nu u=0$, defining what we call the Neumann realization $A_\nu $;
it is selfadjoint with
$D(A_\nu )\subset H^2(\Omega )$. It is well-known
that when $b$ is smooth, then $D(\wA)\subset H^2(\Sigma )$.

\proclaim{Lemma 1.2} When $b\in L_\infty (\Sigma )$, the domain of $\wA$ satisfies 
$$
D(\wA)\subset H^{\frac32}(\Omega )\cap D(\Ama).
$$
\endproclaim

\demo{Proof} When $u\in D(\wA)$, then $u\in H^1(\Omega )$ implies
$\gamma _0u\in H^{\frac12}(\Sigma )\subset L_2(\Sigma
)$. Multiplication by $b$ is continuous on $L_2(\Sigma )$, so $b\gamma
_0u\in L_2(\Sigma )$. Then also $\nu u=b\gamma _0u$ is in $L_2(\Sigma
)$. By the ellipticity of the Neumann problem, $Au\in L_2(\Omega )$
with $\nu u\in L_2(\Sigma )$ imply $u\in H^{\frac32}(\Omega )$. 
\qed\enddemo

When $b$ has some smoothness or piecewise smoothness, we can get more
regularity: It is known that when $b$ is in the Bessel potential space
$H^r_p(\Sigma )$ with $r>(n-1)/p$, $p\ge 2$, then multiplication by $b$
is continuous in $H^s(\Sigma )$ for $|s|\le r$ (cf.\ e.g.\ Johnsen
\cite{J95}). In relation to H\"o{}lder spaces $C^r $ and Besov spaces 
$B^r _{p.q}$ there are inclusions
$$
C^{r+2\delta  }(\Sigma )\hookrightarrow B^{r +\delta 
}_{\infty ,2}(\Sigma )\hookrightarrow B^{r +\delta  }_{p
,2}(\Sigma )\hookrightarrow H^r_p(\Sigma ),\quad \text{ any }\delta  >0,\tag1.10
$$
so also functions in these spaces preserve $H^s(\Sigma  )$ for
$|s|\le r$. (A summary of the relevant facts on function spaces is
given e.g.\ in \cite{AGW10}, Sect.\ 2.) Note that any $\varepsilon >0$
can be included as an $r+2\delta $ by taking $r\in \,]0,\varepsilon
[\,$, $\delta
=(\varepsilon -r)/2$ and $p>(n-1)/r$.

When $X(\Sigma )$ is a function space over $\Sigma $, we
say that $b$ is piecewise in $X$, when
the $(n-1)$-dimensional manifold $\Sigma $ is a union $\Sigma
_1\cup\dots\cup\Sigma _J$ of smooth subsets $\Sigma _j$ with disjoint
interiors (such
that the interfaces are smooth $(n-2)$-dimensional manifolds), and
there are functions $b_j\in X (\Sigma )$,  such that $b$ equals $b_j$ on 
the interior $\Sigma _j^\circ $, for $j=1,\dots,J$.

It is well-known that multiplication by  $1_{\Sigma _j}$ is 
continuous on $H^s(\Sigma )$ for all $|s|<\frac12$.

\proclaim{Proposition 1.3} 

$1^\circ$ Let $b\in H^r_p(\Sigma )$ with $r>(n-1)/p$, $p\ge 2$ (it
holds if $b$ is in one of the spaces in {\rm (1.10)}). Then 
 $D(\wA)\subset
H^{\frac32+r}(\Omega )$ if $r<\frac12$, $D(\wA)\subset H^2(\Omega )$ if $r\ge
\frac12$. 

$2^\circ$ Let  $b$ be piecewise in $H^r_p(\Sigma )$ with 
$r>(n-1)/p$, $p\ge 2$.
Then $D(\wA)\subset
H^{\frac32+r}(\Omega )$ if $r<\frac12$, $D(\wA)\subset
H^{2-\varepsilon }(\Omega )$ for any $\varepsilon >0$ if $r\ge
\frac12$. 

\endproclaim

\demo{Proof} As already noted, $u\in H^1(\Omega )$ implies
$\gamma _0u\in H^{\frac12}(\Sigma )$. In the case $1^\circ$,  multiplication by $b$
preserves $H^s(\Sigma )$ for $|s|\le r$, so 
 $b\gamma
_0u\in H^{\min\{r,\frac12\} }(\Sigma )$. Then also $\nu u=b\gamma _0u$
is in 
$H^{\min\{r,\frac12\} }(\Sigma )$, and now $Au\in L_2(\Omega )$
with $\nu u\in H^{\min\{r,\frac12\} }(\Sigma )$ imply $u\in
H^{\frac32+r}(\Omega )$ if $r<\frac12$, $u\in H^2(\Omega )$ if $r\ge
\frac12$, by the ellipticity of the Neumann problem. 

In the case $2^\circ$,  since $b=\sum_{j=1}^J b_j1_{\Sigma _j}$, multiplication by $b$
maps $H^r(\Sigma )$ into itself if $ r<\frac12$, and into
$H^{\frac12-\varepsilon }$, any $\varepsilon >0$, if $r\ge
\frac12$. Completing the proof as under $1^\circ$, we find that 
$u\in
H^{\frac32+r}(\Omega )$ if $r<\frac12$, $u\in H^{2-\varepsilon }(\Omega )$ if $r\ge
\frac12$. \qed
\enddemo

Let us regard $\wA$ from the point of view of the general extension
theory of \cite{G68}, as recalled in \cite{BGW09}, \cite{G08}, \cite{G10a}. 

We take the Dirichlet realization $A_\gamma $ as the reference
operator, assumed to have a positive lower bound. 
(Seen from the point of view of \cite{G68}, \cite{BL07} 
uses instead  the
Neumann realization $A_\nu $ as the reference operator.) The operator $\wA$ corresponds,
by the general theory of \cite{G68}, to a closed densely defined operator $T\colon
V\to W$, where $V$ and $W$ are closed subsets of $Z=\ker \Ama$, and
$D(T)$ is dense in $V$; and this in turn is carried over by use of the
homeomorphism $\gamma _0:Z\simto H^{-\frac12}(\Sigma )$, to a closed operator $L:X\to Y^*$, with domain
$D(L)$ dense in $ X$, where $X$ and $Y$ are closed subspaces of 
$H^{-\frac12}(\Sigma )$. Here $X=\gamma _0V$, $Y=\gamma _0W$ and $D(L)=\gamma _0D(T)=\gamma _0D(\wA)$.

\proclaim{Proposition 1.4} The operator $L:X\to Y$ corresponding to $\wA$
by \cite{G68} has $X=Y=H^{-\frac12}(\Sigma )$, and acts like
$b-P^0_{\gamma ,\nu }$ with a domain contained in $H^1(\Sigma )$.
When $b$ is real, $L$ is selfadjoint as an unbounded operator from $H^{-\frac12}(\Sigma
)$ to $H^{\frac12}(\Sigma )$.

\endproclaim

\demo{Proof} Besides the description referred to above,
we shall use the observations on operators defined by sesquilinear
forms worked out in  \cite{G70} (and partly recalled in
 \cite{G09},
 Ch.\ 13.2, see in particular Th.\ 13.19). Since the domain of $a_b(u,v)$ equals
 $H^1(\Omega )$, $T$ is defined from a sesquilinear form $t(z,w)$ with domain
 $H^1(\Omega )\cap Z$ dense in $Z$, and hence $V=W=Z$. It follows that
 $X=Y=H^{-\frac12}(\Sigma )$, and $L$ is densely defined and closed as  an
 operator from $H^{-\frac12}(\Sigma )$ to
 $H^{\frac12}(\Sigma )$.
 The adjoint $L^*$ is of the same type and corresponds to $\wA^*$.
When $b$ is real, $\wA$ is selfadjoint as noted above; then
 $L$ is selfadjoint.

In the interpretation of the extension theory, $\wA$ represents the boundary condition
$$
\gamma _0u\in D(L),\quad \Gamma u=L\gamma _0u;
$$
where $\Gamma u=\nu u-P^0_{\gamma ,\nu }\gamma _0u$, so $L\gamma
_0u=\nu u-P^0_{\gamma ,\nu }\gamma _0u$ when $u\in
D(\wA)$. ($P^\lambda _{\gamma ,\nu }$ is the operator mapping
Dirichlet boundary values  to Neumann boundary values for solutions of
$(A-\lambda )u=0$; more on this below.) Since the
functions in $D(\wA)$ also satisfy $\nu u=b\gamma _0u$, we see that 
$L$ acts like 
$$
L\varphi =(b-P^0_{\gamma ,\nu })\varphi .
$$

By Lemma 1.2, $D(\wA)\subset H^{\frac32}(\Sigma )$, so $D(L)=\gamma _0D(\wA)\subset
H^1(\Sigma ) $.
\qed \enddemo

When we replace $A$ by $A-\lambda $, where $\lambda $ is in the
resolvent set  $\varrho (A_\gamma )$ of  $A_\gamma $, we get for the
corresponding operator $L^\lambda $:
$$
L^\lambda \text{ acts like }b-P^\lambda _{\gamma ,\nu },\text{ with
}D(L^\lambda )=D(L)\subset H^1(\Sigma ).
$$

For $\lambda \in \varrho (A_\gamma )\cap \varrho (\wA)$, there holds a
Krein resolvent formula (shown in \cite{BGW09}, Th.\ 3.4):
$$
(\wA-\lambda )^{-1}=(A_\gamma -\lambda )^{-1}+K^\lambda _{\gamma } 
(L^\lambda )^{-1}(K^{\bar\lambda }_{\gamma })^*.\tag1.11
$$
Here $K^\lambda _\gamma $ is the Poisson operator for the Dirichlet
problem, i.e.\ the solution operator $K^\lambda _\gamma \colon
\varphi \mapsto u$ for the problem
$$
(A -\lambda )u=0 \text{ on }\Omega ,\quad \gamma _0u=\varphi \text{
on }\Sigma ;
$$
it maps $H^{s-\frac12}(\Sigma )\to H^{s}(\Omega )$ continuously for all
$s$, and the adjoint $(K^\lambda _\gamma)^*$ maps e.g.\ $L_2(\Omega )$ to $H^{\frac12}(\Sigma )$.

We can use this to show a spectral estimate for $(\wA-\lambda
)^{-1}-(A_\nu  -\lambda )^{-1}$, going via differences with the
Dirichlet resolvent. The argumentation is not the same as that of
\cite{BLLLP10}, which uses a Krein formula based on
the Poisson operator for the Neumann problem.

The spectrum of $A_\gamma $ is contained in a positive halfline
$[c_0,\infty [\,$, and the spectrum of $A_\nu $ is contained in 
a larger halfline $\,]-k,\infty [\,$, cf.\ (1.2). For $\lambda \in
{\Bbb C}\setminus \,]-k,\infty [\, $, the Dirichlet-to-Neumann
operator $P^\lambda _{\gamma ,\nu }=\nu K_\gamma ^\lambda $ is a homeomorphism from
$H^{s}(\Sigma )$ to $H^{s-1}(\Sigma )$ for all $s\in{\Bbb R}$, with
inverse $P^\lambda _{\nu ,\gamma }$, the Neumann-to-Dirichlet
operator. Then we can write
$$
L^\lambda \varphi =(b-P^\lambda _{\gamma ,\nu })\varphi 
=(bP^\lambda _{\nu ,\gamma }-1)P^\lambda _{\gamma ,\nu }\varphi , \text{ for
}\varphi \in D(L).\tag1.12
$$
Since
$P^\lambda _{\nu ,\gamma  }$ is of order $-1$, it is compact in
$L_2(\Sigma )$. Then $bP^\lambda _{\nu ,\gamma }-1$ is a Fredholm
operator in $L_2(\Sigma )$, as noted also in \cite{BLLLP10}.
If $\lambda $ is such that: (1) $L^\lambda $ is invertible (from $D(L)$
to $H^{\frac12}(\Sigma )$), (2) $bP^\lambda _{\nu ,\gamma }-1$ is
invertible in $L_2(\Sigma )$, then the inverse of $L^\lambda $ must
coincide with the inverse of $(bP^\lambda _{\nu ,\gamma }-1)P^\lambda
_{\gamma ,\nu }$ on
$H^\frac12(\Sigma )$.

For $bP^\lambda _{\nu ,\gamma }-1$, we get invertibility as follows:
We have
as a simple application of the principles in \cite{G96} (cf.\ Th.\
2.5.6, (A.25--26)) that
$$
\|P^\lambda _{\gamma ,\nu }\varphi \|_{H^{s,\mu }(\Sigma )}\simeq
\|\varphi \|_{H^{s+1,\mu }(\Sigma )}, \quad \|\varphi \|_{H^{s-1,\mu }(\Sigma )}\simeq \|P^\lambda _{\nu ,\gamma  }\varphi \|_{H^{s,\mu }(\Sigma )},
$$
uniformly in $\mu =|\lambda |^{\frac12}$ for $\lambda \to\infty $ on 
rays in ${\Bbb C}\setminus \rp$; this holds since $P^\lambda _{\gamma
,\nu }$ is parameter-elliptic of order 1 and regularity $+\infty $ on
the rays in ${\Bbb C}\setminus \rp$.
In particular, one has on such a ray $\{\lambda =\mu ^2e^{i\eta }\}$
with $\eta \in \,]0,2\pi [\,$, for $s\in [0,1]$ and $\mu \ge 1$,
$$
\|P^\lambda _{\nu ,\gamma  }\varphi \|_{H^{s}(\Sigma )}+\ang\mu
^s\|P^\lambda _{\nu ,\gamma  }\varphi \|_{L_2(\Sigma )}\le C\min
\{\|\varphi \|_{H^{s-1}(\Sigma )}, \ang\mu ^{s-1}\|\varphi \|_{L_2(\Sigma )}\},
$$
so the norm of $P^\lambda _{\nu ,\gamma  }$ in $L_2(\Sigma )$ is
$O(\ang\mu ^{-1})$ on the ray. Take $\mu _0$ so large that
$\|bP^\lambda _{\nu ,\gamma }\|_{\Cal L(L_2(\Sigma ))}\le \delta <1$
for $\mu \ge \mu _0$, then $bP^\lambda _{\nu
,\gamma }-1$ is invertible as an operator in $L_2(\Sigma )$ for $\mu
\ge \mu _0$, with a bounded inverse  $(bP^\lambda _{\nu,\gamma
}-1)^{-1}$:
$$
(bP^\lambda _{\nu,\gamma
}-1)^{-1}=-1-\sum_{k=1}^\infty (b P^\lambda _{\nu ,\gamma })^k,\text{
converging in }\Cal L(L_2(\Sigma )).\tag1.13
$$
 Then $b-P^\lambda _{\gamma ,\nu } $ has an
inverse 
$$
(b-P^\lambda _{\gamma ,\nu })^{-1}=  P^\lambda _{\nu ,\gamma }(bP^\lambda
_{\nu ,\gamma }-1)^{-1}. \tag1.14
$$

For $L^\lambda $ we know from the extension theory that $L^\lambda $ is
bijective from $D(L)$ to $H^\frac12(\Sigma )$ if and only if $\lambda
\in \varrho (\wA)$. It follows from Lemma 1.1 by a simple geometric
consideration that for each ray
$\{\lambda =\mu ^2e^{i\eta }\}$ with $\eta \in \,]0,2\pi [\,$,
there is a $\mu _1$ such that  such that $\lambda \in \varrho (\wA)$ 
for $\mu \ge \mu _1$. 

For $\mu \ge \max\{\mu _0,\mu _1\}$, both (1) and
(2) are satisfied, so then  
$$
(L^\lambda )^{-1}=(b-P^\lambda _{\gamma ,\nu })^{-1}=  P^\lambda _{\nu ,\gamma }(bP^\lambda
_{\nu ,\gamma }-1)^{-1}\text{ on }H^{\frac12}(\Sigma ). \tag1.15
$$

We note in particular that
$$
D(L^\lambda )=\{\varphi \in H^{1}(\Sigma )\mid (b-P^\lambda _{\gamma ,\nu
})\varphi\in H^{\frac12}(\Sigma ) \},\tag1.16
$$
for such $\lambda $. Now $D(L)=D(L^\lambda )$, and $P^0_{\gamma ,\nu
}-P^\lambda _{\gamma ,\nu }$ is bounded from $H^{-\frac12}(\Sigma )$
to $H^{\frac12}(\Sigma )$ (cf.\ \cite{BGW09}, Rem.\ 3.2), so we
conclude that 
$$
D(L)=\{\varphi \in H^{1}(\Sigma )\mid (b-P^0_{\gamma ,\nu
})\varphi\in H^{\frac12}(\Sigma ) \}.\tag1.17
$$
It follows moreover that (1.16) holds for {\it all }$\lambda \in
\varrho (A_\gamma )$.

This shows the main part of:

\proclaim{Theorem 1.5} The domain of $L$ satisfies {\rm (1.17)}, and
it is also described by  {\rm (1.16)} for any $\lambda \in \varrho
(A_\gamma )$.

On each ray in ${\Bbb C}\setminus \rp$, $\lambda $ is in $\varrho
(\wA)$ and {\rm (1.15)} holds for $|\lambda |$ sufficiently large. For
such $\lambda $,
$$
(\wA-\lambda
)^{-1}-(A_\gamma   -\lambda )^{-1} = K^\lambda _\gamma P_{\nu ,\gamma }(bP^\lambda _{\nu ,\gamma }-1)^{-1}(K^{\bar\lambda }_\gamma )^*.\tag1.18
$$
\endproclaim

\demo{Proof} The statements before formula (1.18) were accounted for
above, and the formula follows by  insertion of (1.15) in {\rm (1.11)}.\qed
\enddemo

\subhead 2. Spectral estimates \endsubhead

Spectral estimates for resolvent differences will now be studied.
A classical reference for the basic concepts is the book of Gohberg
and Krein \cite{GK69}; some particularly
relevant facts were collected in \cite{G84}, supplied with additional
results. We shall include a short summary here:

 For $p>0$, the space $\Cal C_p$ is the Schatten class of compact linear
operators $T$ (in a Hilbert space $H$) with singular value sequences $(s_j(T))_{j\in{\Bbb N}}\in
\ell_p$, and $\frak S_{p}$ denotes the quasi-normed space of
compact operators $T$ with $s_j(T) =O(j^{-1/p})$; here $\frak
S_{p}\subset \Cal C_{p +\varepsilon }$ for all $\varepsilon >0$. 

The rules shown by Ky Fan \cite{F51}$$
s_{j+k-1}(T+T')\le s_j(T)+s_k(T'),\quad s_{j+k-1}(TT')\le
s_j(T)s_k(T'),\tag2.1
$$ 
imply that $\Cal C_p$ and $\frak S_p$ are vector spaces, and that a
product rule holds:$$
\frak S_p\cdot \frak S_q\subset \frak S_{1/(p^{-1}+q^{-1})},\quad \Cal
C_p\cdot \Cal C_q\subset \Cal C_{1/(p^{-1}+q^{-1})}.\tag2.2
$$
Moreover, the rule $$
s_j(ATB)\le \|A\|s_j(T)\|B\|\tag2.3
$$ implies that $\frak
S_p$ and $\Cal C_p$ are preserved under compositions with bounded
operators. We mention two perturbation results:

\proclaim{Lemma 2.1}

$1^\circ$ If $s_j(T)j^{1/p}\to C_0$ and $s_j(T')j^{1/p}\to 0$ for
$j\to\infty$, then $s_j(T+T')j^{1/p}\to C_0$ for $j\to\infty$.

$2^\circ$ If $T=T_M+T'_M$ for each $M\in\Bbb N$, where
$s_j(T_M)j^{1/p}\to C_M$ for $j\to\infty$ and $s_j(T'_M)j^{1/p}\le
\varepsilon _M$ for
$j\in\Bbb N$, with $C_M\to C_0$ and $\varepsilon _M\to 0$ for
$M\to\infty$, then $s_j(T)j^{1/p}\to C_0$ for $j\to\infty$.

\endproclaim

The statement in $1^\circ$ is the Weyl-Ky Fan theorem (cf.\ e.g.\
\cite{GK69} Th.\ II 2.3), and $2^\circ$ is a refinement shown in
\cite{G84}, Lemma 4.2.$2^\circ$.

We also recall that when $\Xi$ is a compact $n'$-dimensional smooth
manifold (possibly with boundary) and $T$ is a bounded linear operator
from $L_2(\Xi)$ to $H^t(\Xi)$ for some $t>0$, then $T\in \frak
S_{n'/t}$ as an operator in $L_2(\Xi)$, with $$
s_j(T)j^{t/n'}\le C\|T\|_{\Cal L(L_2,H^t)},\tag 2.4
$$
$C$ depending only on $\Xi$ and $t$. See \cite{G84}, Lemma
4.4ff.\ for references.

\medskip
The Poisson operator $K^\lambda _\gamma $ is continuous 
from
$H^{s-\frac12}(\Sigma )$ to $H^s(\Omega )$ for all $s\in {\Bbb R}$,
and its adjoint ${K^{\lambda} _\gamma }^*$ is a trace operator of
class 0 and order $-1$ in the pseudodifferential boundary operator 
calculus, hence is
continuous from $H^s(\Omega )$ to $H^{s+\frac12}(\Sigma )$ for
$s>-\frac12$. Then the composition ${K^\lambda _\gamma }^*K^\lambda
_\gamma $ is continuous from $L^2(\Sigma )$ to $H^1(\Sigma )$, so in
view of (2.4),
${K^\lambda _\gamma }^*K^\lambda
_\gamma \in \frak S_{n-1}$ and hence 
$K^\lambda _\gamma \in \frak S_{(n-1)/(1/2)}$, as operators in $L_2(\Sigma)$. 
The singular numbers of
${K^\lambda _\gamma }^*$ have the same behavior.
Moreover, since $P^\lambda _{\nu ,\gamma }$ is a pseudodifferential
operator of order $-1$ on
$\Sigma $, it lies in $\frak S_{n-1}$ when considered as an operator in $L_2(\Sigma)$.

\proclaim{Theorem 2.2} Let $b\in L_\infty(\Sigma )$. For any $ \lambda \in \varrho (\wA)\cap \varrho
(A_\nu )$, 
$$
(\wA-\lambda
)^{-1}-(A_\nu    -\lambda )^{-1}\in \frak S_{(n-1)/3}.\tag 2.5
$$

\endproclaim

\demo{Proof} First assume that $\lambda $ lies so far out on a ray in
$ {\Bbb C}\setminus \rp$ that the statements in Theorem 1.5 are valid.

Applying (1.18) to our $\wA$ and also to the case $b=0$
(the Neumann realization), we find by subtraction:
$$
\aligned
(\wA-\lambda
)^{-1}-(A_\nu    -\lambda )^{-1} &=(\wA-\lambda
)^{-1}-(A_\gamma -\lambda )^{-1} -((A_\nu    -\lambda )^{-1}-(A_\gamma -\lambda )^{-1})\\
&= K^\lambda _\gamma P^\lambda
_{\nu ,\gamma }[(bP^\lambda _{\nu ,\gamma
}-1)^{-1}+1]{K^{\bar\lambda }_\gamma }^*\\
&=K^\lambda _\gamma P^\lambda _{\nu ,\gamma }(bP^\lambda _{\nu
,\gamma }-1)^{-1}bP^\lambda _{\nu ,\gamma }{K^{\bar\lambda }_\gamma }^*.
\endaligned\tag2.6
$$
The last
expression is composed of the operator $K^\lambda _{\gamma }$ in
$\frak S_{(n-1)/(1/2)}$, the adjoint of $K^{\bar\lambda }_\gamma $ with the same property, two factors
$P^\lambda _{\nu ,\gamma }$ in $\frak S_{n-1}$ and the bounded
operators $(bP^\lambda _{\nu
,\gamma }-1)^{-1}$ and $b$, so it belongs to $\frak S_{(n-1)/3}$, by (2.2).

Now let $\lambda '$ be an arbitrary number in $\varrho (\wA)\cap
\varrho (A_\nu )$. We use the following refined resolvent identity as
in \cite{BLLLP10}:
$$
\multline
(S-\lambda ')^{-1}-(T-\lambda ')^{-1}
\\=(1+(\lambda '-\lambda )(T-\lambda ')^{-1})((S-\lambda
)^{-1}-(T-\lambda )^{-1})(1+(\lambda '-\lambda )(S-\lambda ')^{-1}),
\endmultline \tag2.7$$
valid for $\lambda ,\lambda '\in \varrho (T)\cap \varrho
(S)$. Applying it to $S=\wA$ and $T=A_\nu $ for $\lambda $ as above
and $\lambda '\in \varrho (\wA)\cap \varrho (A_\nu )$, we find that 
$(\wA-\lambda ')^{-1}-(A_\nu-\lambda ')^{-1}$ is a composition of an
operator in $\frak S_{(n-1)/3}$ with two bounded operators, hence lies
in $\frak
S_{(n-1)/3}$, as was to be shown.\qed
\enddemo

There is an obvious corollary:

\proclaim{Corollary 2.3}  Let $b_1, b_2\in L_\infty(\Sigma )$, and
denote the corresponding realizations of Robin conditions $\nu
u=b_1\gamma _0u$ resp.\ $\nu u=b_2\gamma _0u$ by $\wA_1$ resp.\ $\wA_2$. For any $ \lambda \in \varrho (\wA_1)\cap \varrho (\wA_2)$, 
$$
(\wA_1-\lambda
)^{-1}-(\wA_2  -\lambda )^{-1}\in \frak S_{(n-1)/3}.\tag 2.8
$$
\endproclaim

\demo{Proof} Write $(\wA_1-\lambda
)^{-1}-(\wA_2  -\lambda )^{-1}$ as the difference between $(\wA_1-\lambda
)^{-1}-(A_\nu   -\lambda )^{-1}$ and $(\wA_2-\lambda
)^{-1}-(A_\gamma   -\lambda )^{-1}$, then the result follows from
Theorem 2.2 (and (2.7)) since $\frak S_p$ is a vector space.\qed 
\enddemo

Formula (1.18) also allows us to show a {\it spectral asymptotics} estimate
for $(\wA-\lambda
)^{-1}-(A_\gamma   -\lambda )^{-1}$ that was obtained in the smooth
case for selfadjoint realizations and negative $\lambda $ in Grubb
\cite{G74}, Sect.\ 8, and Birman and Solomiak \cite{BS80}.
In the former paper it is shown, also for $2m$-order problems, that
the operator is, on the orthogonal complement of its nullspace, {\it isometric}
to an elliptic  pseudodifferential operator on
$\Sigma$ of order $-2m$ (which has the asserted spectral asymptotics); in  the latter
paper exterior domains are included.

\proclaim{Theorem 2.4} Let $b\in L_\infty(\Sigma )$. For any $ \lambda \in \varrho (\wA)\cap \varrho
(A_\gamma )$, 
$$
s_j((\wA-\lambda
)^{-1}-(A_\gamma   -\lambda )^{-1})j^{2/(n-1)}\to C_0^{2/(n-1)} \text{ for }j\to\infty ,\tag2.9
$$
where $C_0$ is the same constant as in the case $b=0$ (where
$\wA=A_\nu$), namely
$$
C_0=
\tfrac1{(n-1)(2\pi )^{n-1}}\int_{\Sigma }\int_{|\xi
'|=1} (\|\tilde k^0\|_{L_2(\rp)}|p^0|^{1/2})^{n-1}
\,d\omega (\xi ') dx';\tag2.10
$$
here $\tilde k^0(x',x_n,\xi ')$ is the principal symbol-kernel of
$K^\lambda _\gamma $ and $p^0(x',\xi ')$ is the principal symbol of $P^\lambda _{\nu ,\gamma }$.

\endproclaim

\demo{Proof} Since the details are perhaps
not very well known, we first give a proof of (2.9)--(2.10) 
in the case $b=0$. We have as an easy special case of (1.18) that
$$
(A_\nu -\lambda
)^{-1}-(A_\gamma   -\lambda )^{-1} = -K^\lambda _\gamma P_{\nu
,\gamma }{K^{\bar\lambda }_\gamma }^*\equiv G_\nu .\tag2.11
$$
This is a singular Green operator with principal boundary symbol operator
$$
g_\nu ^0(x',\xi ',D_n)=-k^0(x',\xi ',D_n)p^0(x',\xi ')k^0(x',\xi ',D_n)^*
$$
in local coordinates, where $k^0$ and $p^0$ are the ($\lambda
$-independent) principal symbols
of $K^\lambda _\gamma $ and $P^\lambda _{\nu ,\gamma }$. At each $(x',\xi ')$,  $k^0(x',\xi ',D_n)\colon {\Bbb C}\to L_2(\rp)$ maps $v\in{\Bbb C}$ to
$\tilde k^0(x',x_n,\xi ')v$, where $\tilde k^0(x',x_n,\xi ')\in \Cal
S(\crp)$ is the symbol-kernel. In the case $A=-\Delta $ it equals
$e^{-|\xi '|x_n}$, and it has a similar structure for general $A$
(cf.\ e.g.\ \cite{GS01}, Sect.\ 2.d). The
operator $k^0(x',\xi ',D_n)^*\colon L_2(\rp)\to {\Bbb C}$ maps
$u(x_n)$ to $(u,\tilde k^0)_{L_2(\rp)}$. Thus $k^0(x',\xi ',D_n)^*k^0(x',\xi ',D_n)
$ is the multiplication by $\|\tilde k^0\|^2_{L_2(\rp)}$, and $k^0(x',\xi ',D_n)k^0(x',\xi ',D_n)^*
$ is the rank 1 operator mapping $u$ to $(u,\tilde k^0)\tilde
k^0$. The latter operator has the sole eigenvector $\tilde k^0_1=\tilde
k^0/\|\tilde k^0\|$ with a positive eigenvalue $\|\tilde k^0\|^2$ (besides
eigenvectors in the nullspace), so its trace
equals the eigenvalue. The middle factor $p^0$ is just multiplication
by a scalar; for $A=-\Delta $, it equals $-|\xi '|^{-1}$.

By \cite{G84},
Th.\ 4.10, since $G_\nu $ is a singular Green operator of order $-2$
and class 0, 
$$
s_j(G_\nu )j^{2/(n-1)}\to C(g^0_ \nu )^{2/(n-1)}\text{ for }j\to
\infty ,\tag2.12 
$$
where 
$$
C(g^0_\nu )=\tfrac1{(n-1)(2\pi )^{n-1}}\int_{\Sigma }\int_{|\xi
'|=1}\tr\big[\big(g^0_\nu (x',\xi ',D_n)^*g^0_\nu (x',\xi
',D_n)\big)^{(n-1)/4}\big]\,d\omega (\xi ') dx'.\tag2.13
$$
Here
$$
\aligned
g^0_\nu (x',\xi ',&D_n)^*g^0_\nu (x',\xi
',D_n)\\
&=k^0(x',\xi ',D_n)\bar p^0(x',\xi ')k^0(x',\xi
',D_n)^*k^0(x',\xi ',D_n)p^0(x',\xi ')k^0(x',\xi ',D_n)^*\\
&=\|\tilde k^0(x',x_n,\xi
')\|_{L_2(\rp)}^2| p^0(x',\xi ')|^2k^0(x',\xi ',D_n)k^0(x',\xi ',D_n)^*.
\endaligned
$$
This is a rank 1 operator with eigenvalue $\|\tilde
k^0\|^4_{L_2}|p^0|^2$, so 
$$
\tr\big[\big(g^0_\nu (x',\xi ',D_n)^*g^0_\nu (x',\xi
',D_n)\big)^{(n-1)/4}\big]=(\|\tilde k^0\|_{L_2}|p^0|^{1/2})^{n-1},
$$
and (2.10) follows.

Now the case of general $b$:
For large $\lambda $ on rays in ${\Bbb C}\setminus \rp$ as in Theorem 1.5 we write formula (1.13) as
$$
(bP^\lambda _{\nu,\gamma
}-1)^{-1}=-1-b P^\lambda _{\nu ,\gamma }S,\text{
where }S=\sum_{k=0}^\infty (b P^\lambda _{\nu ,\gamma })^k \in\Cal L(L_2(\Sigma )).\tag2.14
$$
Then we have from (1.18):
$$
\aligned
(\wA-\lambda
)^{-1}-(A_\gamma   -\lambda )^{-1} &= K^\lambda _\gamma P_{\nu ,\gamma
}(-1-bP^\lambda _{\nu ,\gamma } S){K^{\bar\lambda }_\gamma
}^*\\
&= -K^\lambda _\gamma P_{\nu ,\gamma
}{K^{\bar\lambda }_\gamma
}^*-
 K^\lambda _\gamma P_{\nu ,\gamma
}bP^\lambda _{\nu ,\gamma } S{K^{\bar\lambda }_\gamma }^*.
\endaligned\tag2.15
$$
The first term equals $(A_\nu -\lambda
)^{-1}-(A_\gamma   -\lambda )^{-1} $ and satisfies the
spectral asymptotics estimate (2.9) with (2.10). The second term is in $\frak
S_{(n-1)/3}$, in view of the mapping properties of its factors, as in
the proof of Theorem 2.2. By Lemma 2.1.$1^\circ$, it follows that the sum of the two
terms has the asymptotic behavior (2.9).

General $\lambda \in \varrho (\wA)\cap \varrho (A_\gamma )$ are
included by use of the
resolvent identity (2.7), which gives the operator as a sum of a
term with the behavior (2.9) and terms in $\frak
S_{(n-1)/(2+t)}$ with $t>0$,
using that $(A_\gamma-\lambda)^{-1}\in \frak S_{n/2}$ and
$(\wA-\lambda)^{-1}\in \frak S_{n/(3/2)}$. Then Lemma 2.1.$1^\circ$
 applies to show (2.9) for the sum.
\qed
\enddemo

Spectral asymptotics estimates for the resolvent difference (2.5) are
harder to get at, since $b$ here enters in the principal part of the
operator. However, with a little smoothness of $b$ we can obtain the
spectral estimate by reduction to a case that allows an approximation 
procedure.

We consider the resolvent difference of two general Robin problems from the
start, since the asymptotic property is not in general additive.

\proclaim{Theorem 2.5} Assume that $b_1,b_2\in H^r_p(\Sigma )$, where
 $r>0$ and  
 $p>(n-1)/r$, $p\ge 2$; this
holds if the $b_i$ are in one of the spaces in {\rm (1.10)}. Define $\wA_i$ as in Corollary
{\rm 2.3}. Then for $\lambda \in \varrho
 (\wA_1)\cap \varrho (\wA_2 )$, 
$$
s_j((\wA_1-\lambda
)^{-1}-(\wA_2  -\lambda )^{-1} 
)j^{3/(n-1)}\to C(g^0)^{3/(n-1)}\text{ for }j\to \infty ,\tag2.16
$$
where  
$$
C(g^0)=\tfrac1{(n-1)(2\pi )^{n-1}}\int_{\Sigma }\int_{|\xi
'|=1}(\|\tilde k^0\|^2_{L_2(\rp)}|p^0|^2|b_2-b_1|)^{(n-1)/3}
\,d\omega (\xi ') dx'.\tag2.17
$$
\endproclaim

\demo{Proof} 
First let $\lambda $ be large on a ray in ${\Bbb C}\setminus\rp$ such
that Theorem 1.5 applies to $\wA_1$ and $\wA_2$. 
Using (2.14) in the form
$$
(b_iP^\lambda _{\nu,\gamma
}-1)^{-1}=-1-b_iP^\lambda _{\nu ,\gamma }-(b _iP^\lambda _{\nu ,\gamma })^2S_i
$$
 we have that
$$
(b_1P^\lambda _{\nu,\gamma
}-1)^{-1} - (b_2 P^\lambda _{\nu,\gamma
}-1)^{-1} 
=(b_2-b_1)P^\lambda _{\nu ,\gamma }-(b _1P^\lambda _{\nu ,\gamma })^2S_1
-(b _2P^\lambda _{\nu ,\gamma })^2S_2.
$$
Then we get using (2.6):
$$\aligned
(\wA_1&-\lambda
)^{-1}-(\wA_2    -\lambda )^{-1}=(\wA_1-\lambda
)^{-1}-(A_\nu   -\lambda )^{-1}-((\wA_2-\lambda
)^{-1}-(A_\nu   -\lambda )^{-1})\\ 
&= K^\lambda _\gamma P^\lambda
_{\nu ,\gamma }[(b_1P^\lambda _{\nu ,\gamma
}-1)^{-1}+1]{K^{\bar\lambda }_\gamma }^*-
K^\lambda _\gamma P^\lambda
_{\nu ,\gamma }[(b_2P^\lambda _{\nu ,\gamma
}-1)^{-1}+1]{K^{\bar\lambda }_\gamma }^*\\
&=K^\lambda _\gamma P^\lambda _{\nu ,\gamma }(b_2-b_1)P^\lambda _{\nu
,\gamma }{K^{\bar\lambda }_\gamma }^*-K^\lambda _\gamma P^\lambda
_{\nu ,\gamma }(b_1P^\lambda _{\nu ,\gamma })^2
S_1{K^{\bar\lambda
}_\gamma }^* +K^\lambda _\gamma P^\lambda
_{\nu ,\gamma }(b_2P^\lambda _{\nu ,\gamma })^2
S_2{K^{\bar\lambda
}_\gamma }^*\\
&=G+F_1+F_2.
\endaligned\tag2.18
$$
In the terms $F_i$ we use for one of the factors $b_iP^\lambda_{\nu ,\gamma }$ that $b_i$ preserves $H^s(\Sigma )$ for
$|s|\le r$ (see the text before Proposition 1.3), so that  
$b_iP^\lambda_{\nu ,\gamma }$ 
maps $L_2(\Sigma )$ continuously into $H^{r'}(\Sigma
)$, $r'=\min\{r,1\}$. So this factor is in $\frak
S_{(n-1)/r'}$, together with the usual two factors in $\frak S_{(n-1)/(1/2)}$ and
two factors in $\frak S_{n-1}$, whereby the full composed operator $F_i$ is
in $\frak S_{(n-1)/(3+r')}$. It will not influence the spectral
asymptotics.

In the term $G$, let us denote $b_2-b_1=b$. We write $b$ for each $M\in {\Bbb N}$ as a sum
$$
b=b_M+b'_M,\tag2.19
$$
where $b_M\in C^\infty (\Sigma )$ and $\sup_{x'\in\Sigma }|b'_M(x')|\le 1/M$; this is
possible since $b$ is continuous on the smooth compact manifold
$\Sigma $. Accordingly, we write $G=G_M+G'_M$ with
$$
G_M=K^\lambda _\gamma P^\lambda _{\nu ,\gamma }b_MP^\lambda _{\nu
,\gamma }{K^{\bar\lambda }_\gamma }^*,\quad G'_M=K^\lambda _\gamma P^\lambda _{\nu ,\gamma }b'_MP^\lambda _{\nu
,\gamma }{K^{\bar\lambda }_\gamma }^*.
$$
 
Here $G'_M$ is a composition of fixed operators with the usual
$\frak S_{p}$-properties and a factor $b'_M$ whose norm in $\Cal
L(L_2(\Sigma ))$ is $\le 1/M$; this implies that 
$$
\sup_js_j(G'_M)j^{3/(n-1)}\le C/M, \text{ all }M,\tag2.20
$$
for a suitable constant $C$, in view of (2.3).

The term $G_M$ is treated by application of the
tools in \cite{G84}. Since $b_M\in C^\infty $, $G_M$ is a genuine singular
Green operator of order $-3$ and class 0, with polyhomogeneous
symbol. The principal symbol $g^0_M$ is the symbol of the boundary
symbol operator (in local coordinates)
$$
g^0_M(x',\xi ',D_n)=k^0(x',\xi ',D_n)p^0(x',\xi ')b_M(x')p^0(x',\xi
')k^0(x',\xi ',D_n)^*.\tag2.21
$$
It follows from \cite{G84},
Th.\ 4.10, that
$$
s_j(G_M)j^{3/(n-1)}\to C(g^0_M)^{3/(n-1)}\text{ for }j\to \infty ,\tag2.22
$$
where 
$$
C(g^0_M)=\tfrac1{(n-1)(2\pi )^{n-1}}\int_{\Sigma }\int_{|\xi
'|=1}\tr\big[\big(g^0_M(x',\xi ',D_n)^*g^0_M(x',\xi
',D_n)\big)^{(n-1)/6}\big]\,d\omega (\xi ') dx'.\tag2.23
$$
As in the analysis of ${g^0_\nu}^*g^0_\nu  $ in the proof of Theorem
2.4, now with
the middle factor $p^0$ replaced by $p^0b_Mp^0$, we find that
$$
\aligned
\tr\big[\big(g^0_M (x',\xi ',D_n)^*g^0_M (x',\xi
',D_n)\big)^{(n-1)/6}\big]&=(\|\tilde
k^0\|^4_{L_2}|p^0|^4|b_M|^2)^{(n-1)/6}\\
&=(\|\tilde k^0\|^2_{L_2}|p^0|^2|b_M|)^{(n-1)/3},
\endaligned
$$
and hence$$
C(g^0_M)=
\tfrac1{(n-1)(2\pi )^{n-1}}\int_{\Sigma }\int_{|\xi
'|=1} (\|\tilde k^0\|^2_{L_2(\rp)}|p^0|^2|b_M|)^{(n-1)/3}
\,d\omega (\xi ') dx'.\tag2.24
$$
When $M\to\infty $,
$b_M(x')\to b(x')$ uniformly in $x'$, so   
$$
\multline
C(g^0_M)\to C(g^0),\text{ where } \\
C(g^0)=\tfrac1{(n-1)(2\pi )^{n-1}}\int_{\Sigma }\int_{|\xi
'|=1}(\|\tilde k^0\|^2_{L_2(\rp)}|p^0|^2|b|)^{(n-1)/3}
\,d\omega (\xi ') dx',\\
\text{ with }b=b_2-b_1.
\endmultline\tag2.25
$$

Now we first apply Lemma 2.1.$2^\circ$ to the decompositions
$G=G_M+G'_M$; this shows that $G$ has the spectral behavior in
(2.16). When $F_1$ and $F_2$ are added to $G$, we can use Lemma
2.1.$1^\circ$ 
to conclude that also $G+F_1+F_2$ has
the spectral behavior in (2.16).

Finally, general $\lambda \in \varrho (\wA)\cap \varrho (A_\nu )$ are
included by use of the resolvent formula (2.7) as in the preceding
proof.
\qed

\enddemo

In the case $A=-\Delta $, where $\tilde k^0$ and $p^0$ are independent
of $x'$,  the formula for $C(g^0)$ reduces to a constant times $\int_{\Sigma
}|b_2-b_1|^{(n-1)/3}\, dx'$.

\subhead 3. Coefficients with jumps \endsubhead

It possible to extend the result of Theorem 2.5 to cases where $b$ has
jump discontinuities, by use of special results for pseudodifferential
operators (from here on abbreviated to $\psi $do's). In showing this,
we also supply the general knowledge on spectral asymptotics for
$\psi $do's multiplied with nonsmooth functions.

Let $\Xi $ be a compact $n'$-dimensional $C^\infty $-manifold without
boundary, and assume that it is divided by a smooth
$(n'-1)$-dimensional hypersurface into two subsets $\Xi _+$ and $\Xi
_-$ ($n'$-dimensional $C^\infty $-manifolds with boundary) such that 
$\Xi =\Xi _+\cup\Xi _-$, $\Xi _+^\circ\cap \Xi _-^\circ=\emptyset$,
$\partial\Xi _+=\partial\Xi _-$. (Since the sets need not be connected,
this covers the situation of $J$ smooth subsets described before
Proposition 1.3.) We denote by $r^\pm$ the restrictions from $\Xi
$ to $\Xi _\pm$, and by $e^\pm$ the extension-by-zero operators from
functions on $\Xi _\pm$ to functions on $\Xi $:
$$
e^\pm u=\cases u\text{ on }\Xi _\pm\\ 0\text{ on }\Xi _\mp.\endcases
$$
Multiplication by the characteristic function $1_{\Xi _+}$ for $\Xi _+$
can also be written $e^+r^+$; similarly $1_{\Xi _-}=e^-r^-$.
 
It is well-known (as recalled e.g.\ in \cite{G84}, Lemma 4.5) that
when $P$ is an $N\times N$-matrix formed  classical $\psi $do on $\Xi $ of negative order $-t$, then 
it satisfies the spectral asymptotics formulas for $j\to\infty$:
$$
\aligned
s_j(P)j^{t/n'}& \to C(p^0)^{t/n'} \text{ in general},\\
\pm\lambda ^{\pm}_j(P)j^{t/n'}& \to C^{\pm}(p^0)^{t/n'}\text{ if $P$
is selfadjoint} ,
\endaligned\tag3.1
$$
where, respectively, 
$$
\aligned
C(p^0) &=
\tfrac1{n'(2\pi )^{n'}}\int_{\Xi }\int_{|\xi
|=1}\tr\big[
\big(p^{0}(x,\xi )^*p^0(x,\xi )\big)^{n'/2t}
\big]\,d\omega (\xi ) dx,\\
C^{\pm}(p^0) &=
\tfrac1{n'(2\pi )^{n'}}\int_{\Xi }\int_{|\xi
|=1}{\sum}_{\operatorname{ev.}\gtrless 0 }\big(
\pm\lambda _j^{\pm}(p^{0}(x,\xi ))^{n'/t}
\big)\,d\omega (\xi ) dx.
\endaligned\tag3.2
$$

Let us also recall the result of Laptev \cite{L77, L81}:

\proclaim{Proposition 3.1} Let $P$ be a classical pseudodifferential
operator on $\Xi $ of negative order $-t$. Then $1_{\Xi _+}P1_{\Xi _-}\in \frak
S_{(n'-1)/t}$.
\endproclaim

(Expressed in local coordinates, this means that the operator whose
kernel is the restriction of the kernel of $P$ to the second or fourth
quadrant, picks up the boundary dimension in its spectral behavior. For
$\psi $do's having the transmission property at $\partial\Xi _+$, this
is confirmed by the results of \cite{G84}.)

The rules in the following are valid also for $N\times N$-matrix
formed operators $P$ and factors $b$, and would then need a trace
indication tr in the integrals; we leave this aspect out here for simplicity.

\proclaim{Theorem 3.2} Let $P$ be a classical pseudodifferential
operator of negative order $-t$, such that $(Pu,u)\ge 0$ for $u\in L_2(\Xi )$. Then $P_{(+)}=1_{\Xi _+}P1_{\Xi _+}$
satisfies the spectral asymptotics formula
$$
s_j(P_{(+)})j^{t/n'} \to c(P_{(+)})^{t/n'}\text{ for }j\to\infty ,\tag3.3
$$
where 
$$
\aligned
c(P_{(+)}) &=
\tfrac1{n'(2\pi )^{n'}}\int_{\Xi _+}\int_{|\xi
|=1}
\big(p^{0}(x,\xi )^*p^0(x,\xi )\big)^{n'/2t}
\,d\omega (\xi ) dx\\
&=
\tfrac1{n'(2\pi )^{n'}}\int_{\Xi _+}\int_{|\xi
|=1}
p^0(x,\xi )^{n'/t}
\,d\omega (\xi ) dx.
\endaligned
\tag3.4
$$

\endproclaim

\demo{Proof} The principal symbol $p^0$ is $\ge 0$; which explains the
second identity in (3.4).
Introduce two $C^\infty $ cutoff functions $\zeta
_1$ and $\zeta _2$
taking values in $[0,1]$ such that $\zeta _1=1 $ on $\Xi _+$ and
vanishes outside a neighborhood of $\Xi _+$, and $\zeta _2=0$ on $\Xi _-$
and is 1 outside a neighborhood of $\Xi _-$. We shall then compare
$P_{(+)}$ with the operators (all are compact in $L_2(\Xi )$)
$$
P_1=\zeta _1P\zeta _1\text{ and } P_2=\zeta _2P\zeta _2.
$$

When $u\in L_2(\Xi )$, denote $e^\pm r^\pm u=u_\pm$. We have for
$P_1$, since $\zeta _1u_+=u_+$:
$$
\aligned
(P_1u,u)&=(P_1u_+,u_+)+(P_1u_+,u_-)+(P_1u_-,u_+)+(P_1u_-,u_-)\\
&=(P_{(+)}u,u)+(Ru,u)+(P\zeta _1u_-,\zeta _1u_-),
\endaligned
$$
where $R=1_{\Xi _-}P_1I_{\Xi _+}+1_{\Xi _+}P_1I_{\Xi _-}$. Since
$P_1$ is a classical $\psi $do of order $-t$ on $\Xi $, it has the
spectral behavior in (3.1)--(3.2) with the limit $C(p_1^0)^{t/n'}$;
here 
$$
C(p^0_1)=
\tfrac1{n'(2\pi )^{n'}}\int_{\supp\zeta _1}\int_{|\xi
|=1}
(\zeta _1p^0(x,\xi )\zeta _1)^{n'/t}
\,d\omega (\xi ) dx.
$$
Moreover,
$R$ is of the type considered in Proposition 3.1, hence lies in $\frak
S_{(n'-1)/t}$. Then by Lemma 2.1.$1^\circ$, $P_1-R$ likewise has the
spectral behavior in (3.1)--(3.2) with the limit $C(p_1^0)^{t/n'}$. Now
observe that since  $P$ is nonnegative, $(P\zeta _1u_-,\zeta _1u_-)\ge
0$ for all $u\in L_2(\Xi )$. Thus we have:
$$
(P_{(+)}u,u)\le ((P_1-R)u,u),\text{ for all }u\in L_2(\Xi ).\tag3.5
$$
Both operators $P_{(+)}$ and $P_1-R$ are selfadjoint nonnegative, so the
$s$-numbers are the same as the eigenvalues, and the minimum-maximum
principle implies in view of (3.5) that 
$$
s_j(P_{(+)})\le s_j(P_1-R),\text{ for all }j.\tag3.6
$$
It then follows from the limit property of the $s_j(P_1-R)$ that 
$$
{\lim\sup}_{j\to\infty }s_j(P_{(+)})j^{t/n'}\le C(p^0_1)^{t/n'}.\tag3.7
$$

For the comparison with $P_2$ we write, using that $\zeta _2u_+=\zeta _2u$,
$$\aligned
(P_{(+)}u,u)&=(\zeta _2P\zeta _2u_+,u_+)+((1-\zeta _2)P(1-\zeta
_2)u_+,u_+)
+((1-\zeta _2)P\zeta _2u_+,u_+)\\
&\qquad\qquad+(\zeta _2P(1-\zeta _2)u_+,u_+)\\
&\ge(\zeta _2P\zeta _2u,u)
+((1-\zeta _2)P\zeta _2u,u_+)+(\zeta _2P(1-\zeta _2)u_+,u)\\
&=((\zeta _2P\zeta _2+(1-\zeta _2)P\zeta _2+\zeta _2P(1-\zeta _2))u,u)+(R_1u,u),
\endaligned$$
where $R_1$ is a sum of terms as in Proposition 3.1. Then since
$s_j(P_{(+)})=\lambda _j(P_{(+)})\ge \lambda ^+_j(\zeta _2P\zeta _2+(1-\zeta _2)P\zeta _2+\zeta _2P(1-\zeta _2)+R_1)
$, 
 $$
{\lim\inf}_{j\to\infty }s_j(P_{(+)})j^{t/n'}\ge C^+(\zeta _2p^0\zeta _2+(1-\zeta _2)p^0\zeta _2+\zeta _2p^0(1-\zeta _2))^{t/n'}.\tag3.8
$$

Since $C(p^0_1)$ and $C^+(\zeta _2p^0\zeta _2+(1-\zeta _2)p^0\zeta _2+\zeta _2p^0(1-\zeta _2))$ come arbitrarily close to $c(P_{(+)})$
when the support of $\zeta _1$ shrinks towards $\Xi _+$ and 
the support of $1-\zeta _2$ shrinks towards $\Xi _-$, we conclude that
(3.3) with (3.4) holds.\qed 
\enddemo

This leads to a result on compositions of $\psi $do's with
discontinuous factors, which seems to have an interest
in itself:

\proclaim{Theorem 3.3} Let $P$ be an operator composed of ${l}$
classical pseudodifferential operators $P_1,\dots, P_{l}$ of negative orders
$-t_1,\dots, -t_{l}$ and  ${l}+1$ functions $b_1,\dots,b_{{l}+1}$ that are
piecewise continuous on $\Xi $ with possible jumps at $\partial\Xi _+$
(so the $b_k$ extend to continous funcxtions on $\Xi _+$ and on $\Xi _-$);
$$
P=b_1P_1\dots b_{l}P_{l}b_{l+1}.\tag3.9
$$
Let $t=t_1+\dots+t_{l}$. Then $P$ has the spectral behavior:
$$
s_j(P)j^{t/n'} \to c(P)^{t/n'}\text{ for }j\to\infty ,\tag3.10
$$
where 
$$
\aligned
c(P) &=
\tfrac1{n'(2\pi )^{n'}}\int_{\Xi }\int_{|\xi
|=1}
\big(\bar b_{{l}+1}(x)p_{l}^{0}(x,\xi )^*\dots p_1^{0}(x,\xi )^*\bar
b_1(x)\cdot\\
&\qquad\qquad \qquad\cdot
b_1(x)p_1^0(x,\xi )\dots p_{l}^0(x,\xi )b_{l+1}(x)
\big)^{n'/2t}
\,d\omega (\xi ) dx\\
&=\tfrac1{n'(2\pi )^{n'}}\int_{\Xi }\int_{|\xi
|=1}
|b_1\dots b_{l+1}p^0_1\dots p^0_l|
^{n'/t}
\,d\omega (\xi ) dx.
\endaligned\tag3.11
$$
\endproclaim

\demo{Proof} We can write
$$
P^*P=\bar b_{{l}+1}P_{l}^*\dots P_1^*\bar b_1 
b_1P_1\dots P_{l}b_{l}
=1_{\Xi _+}P^*P1_{\Xi _+}+1_{\Xi
_-}P^*P1_{\Xi _-}+R,
$$
where $R=1_{\Xi _+}\bar b_{{l}+1}P_{l}^*\dots P_{l}b_{l}1_{\Xi
_-}+1_{\Xi _-}\bar b_{{l}+1}P_{l}^*\dots P_{l}b_{l}1_{\Xi
_+}$. Inserting $1=1_{\Xi _+}+1_{\Xi _-}$ at each factor $b_k$ or
$\bar b_k$ in $R$ and multiplying out, we obtain it as a sum of terms of order $-t$, each containing at least one
factor of the type in Proposition 3.1. Thus $R\in \frak
S_{n'/(t+\delta )}$ with a $\delta >0$. For the term $1_{\Xi
_+}P^*P1_{\Xi _+}$,
 we proceed as in Theorem 2.5. We can assume that $b_k$ is extended
 from $\Xi _+$ to a continuous function $b_k$ on $\Xi $. Each $b_k$ is
approximated by a uniformly convergent sequence $b_{kM}$ of $C^\infty
$-functions on $\Xi $. For each $M$, 
$$
P_M^*P_M=\bar b_{{l}+1,M}P_{l}^*\dots P^*_{1}\bar b_{1M}b_{1M}P_{1}\dots b_{{l}M}P_{{l}}b_{{l}+1,M}
$$ 
is a classical nonnegative $\psi $do of order $-t$, so Theorem 3.2 applies 
to the operator with $1_{\Xi _+}$ before and after,
and gives the corresponding spectral asymptotics
formula. Since $P_M^*P_M-P^*P$ can be written as a sum of terms where
each has a small factor $b_{kM}-b_k$ or $\bar b_{kM}-\bar b_k$, we have for
 $M\to\infty $ that 
$$
\sup_js_j(1_{\Xi _+}P_M^*P_M1_{\Xi _+}-1_{\Xi _+}P^*P1_{\Xi_+})j^{t/n'}\to 0.\tag3.12
$$
Then Lemma 2.1.$2^\circ$ implies a spectral asymptotics formula for
$1_{\Xi _+}P^*P1_{\Xi_+}$, with the constant as in (3.11) but integrated
over $\Xi _+$. --- There is a similar result for $1_{\Xi
_-}P^*P1_{\Xi_-}$, relative to $\Xi _-$. 

Now since $L_2(\Xi )$ identifies with the
orthogonal sum of $L_2(\Xi _+)$ and $L_2(\Xi _-)$, the spectra are
simply superposed when the operators are added together. The statement
$\lambda _j(T)j^{t/n'}\to c(T)^{t/n'}$ for $j\to\infty $ is equivalent 
with $N'(a;T)a^{n'/t}\to c(T)$ for $a\to\infty $, 
where $N'(a;T)$ 
counts the number of
eigenvalues in $[1/a,\infty [\,$;  superposition of the
spectra means addition of the counting functions. (More on
counting functions e.g.\ in \cite{G96}, Sect.\ A.6.) 
Thus $1_{\Xi _+}P^*P1_{\Xi_+}+1_{\Xi_-}P^*P1_{\Xi_-}$ has a spectral asymptotics behavior where the constant
is obtained  by adding the integrals for $1_{\Xi _+}P^*P1_{\Xi_+}$ and
$1_{\Xi_-}P^*P1_{\Xi_-}$, so it is as described in
(3.9)--(3.11). By
Lemma 2.1.$1^\circ$, the behavior keeps this form when we add $R$ to
the operator. 
\qed
\enddemo

A similar theorem holds for matrix formed operators $P_k$ and factors
$b_k$, with $c(P)$ defined by the first expression in (3.11); here of
course it cannot be reduced to the second expression unless all the
factors commute.

A special case of the situation in Theorem 3.3 is the case
of $bP$, where $P$ is a classical $\psi $do and $b$ is a piecewise
continuous function. We need a case with interspersed factors $b_k$ in
our application below. 

We can now show:

\proclaim{Theorem 3.4} The conclusion of Theorem {\rm 2.5} holds also
when $b_1$ and $b_2$ are piecewise in $H^r_p(\Sigma )$ for some $r>0$
as in Theorem {\rm 2.5},
$b_2-b_1$ having jumps at a smooth hypersurface.
\endproclaim

\demo{Proof} We use again the decomposition in (2.18):
$$\aligned
&(\wA_1-\lambda
)^{-1}-(\wA_2    -\lambda )^{-1}=G+F_1+F_2,\text{ with }
G=K^\lambda _\gamma P^\lambda _{\nu ,\gamma }(b_2-b_1)P^\lambda _{\nu
,\gamma }{K^{\bar\lambda }_\gamma }^*,\\ 
&F_1=-K^\lambda _\gamma P^\lambda
_{\nu ,\gamma }(b_1P^\lambda _{\nu ,\gamma })^2
S_1{K^{\bar\lambda
}_\gamma }^* ,\quad F_2=K^\lambda _\gamma P^\lambda
_{\nu ,\gamma }(b_2P^\lambda _{\nu ,\gamma })^2
S_2{K^{\bar\lambda
}_\gamma }^*,
\endaligned
$$
and $F_1$ and $F_2$ are handled as after (2.18), using that
$b_iP^\lambda _{\nu ,\gamma }$ maps $L_2(\Sigma )$ into $H^{r'}(\Sigma )$,
$r'=\min\{r, \frac12-\varepsilon \}$. Then they are
in $\frak S_{(n-1)/(3+r')}$. We denote again $b_2-b_1=b$.

For $G$ we proceed as follows: Let $\lambda $ be large negative, so
that Theorem 1.5 holds. Since $\lambda $ is real, $K^{\bar\lambda }_\gamma =K^\lambda
_\gamma $, and $P^\lambda _{\nu ,\gamma }$ is selfadjoint. The $j$-th eigenvalue of $G^*G$ satisfies
$$
\lambda _j(G^*G)=\lambda _j(K^\lambda _\gamma P^\lambda _{\nu ,\gamma
}\bar bP^\lambda _{\nu
,\gamma }{K^{\lambda }_\gamma }^*K^\lambda _\gamma P^\lambda _{\nu ,\gamma }bP^\lambda _{\nu
,\gamma }{K^{\lambda }_\gamma }^*).
$$
Here ${K^{\lambda }_\gamma }^*K^\lambda _\gamma$ equals a selfadjoint
$\psi $do $P_1$ of order $-1$; it is nonnegative on $L_2(\Sigma )$ and
injective, since $K^\lambda _\gamma$ is injective:
$$
(P_1\varphi ,\varphi )_{L_2(\Sigma )}=({K^{\lambda }_\gamma }^*K^\lambda _\gamma
\varphi ,\varphi )_{L_2(\Sigma )}=\|K^\lambda _\gamma \varphi \|^2_{L_2(\Omega )}\ge c\|\varphi \|^2_{H^{-\frac12}(\Sigma )},
$$
hence elliptic. It follows from Seeley \cite{S67} that $P_1$ has a
squareroot $P_2=P_1^\frac12$ which is a classical elliptic $\psi $do of order
$-\frac12$. Then we find, applying the general formula
$$\lambda _j(TT')=\lambda _j(T'T),\tag3.13
$$
with $T=K^\lambda _\gamma P^\lambda _{\nu ,\gamma
}\bar bP^\lambda _{\nu
,\gamma }P_2$, $T'=P_2 P^\lambda _{\nu ,\gamma }bP^\lambda _{\nu
,\gamma }{K^{\lambda }_\gamma }^*$, that
$$
\aligned
\lambda _j(G^*G)&=\lambda _j(K^\lambda _\gamma P^\lambda _{\nu ,\gamma
}\bar bP^\lambda _{\nu
,\gamma }P_2P_2 P^\lambda _{\nu ,\gamma }bP^\lambda _{\nu
,\gamma }{K^{\lambda }_\gamma }^*)\\
&=\lambda _j(P_2 P^\lambda _{\nu ,\gamma }bP^\lambda _{\nu
,\gamma }{K^{\lambda }_\gamma }^*K^\lambda _\gamma P^\lambda _{\nu ,\gamma
}\bar bP^\lambda _{\nu
,\gamma }P_2
)\\
&=\lambda _j(P_2 P^\lambda _{\nu ,\gamma }bP^\lambda _{\nu
,\gamma }P_1 P^\lambda _{\nu ,\gamma
}\bar bP^\lambda _{\nu
,\gamma }P_2
).\endaligned
$$
The operator $Q=P_2 P^\lambda _{\nu ,\gamma }bP^\lambda _{\nu
,\gamma }P_1 P^\lambda _{\nu ,\gamma
}\bar bP^\lambda _{\nu
,\gamma }P_2$ is an operator to which Theorem 3.4 applies, and it
gives a spectral asymptotics formula with the constant defined as in
(3.11), with $n'=n-1$. Since $p_1^0=\|\tilde k^0\|^2_{L_2}$,
$p_2^0=\|\tilde k^0\|_{L_2}$, the formula can be
rewritten in the form (2.17).

The proof is now completed in the same way as in the proof of Theorem 2.5.\qed  
\enddemo

The results can be extended to exterior domains by the method of \cite{G10}.

In a forthcoming paper we shall treat the question of spectral
asymptotics for the {\it mixed problem} for $-\Delta +a_0$, where the boundary
condition jumps from a Dirichlet condition to a Neumann condition at a
smooth hypersurface of $\Sigma $. Here we moreover need to draw on the
analyses of nonstandard pseudodifferential operators, as in Shamir
\cite{S63}, Eskin \cite{E81}, Birman and Solomiak \cite{BS77} 
and many later works.


\Refs
\widestnumber\key{[BLLLP10]}

\ref\no[AGW10]\by H. Abels, G. Grubb and I. Wood \paper Extension theory and  Kre\u\i{}n-type resolvent
  formulas for nonsmooth boundary value
  problems \finalinfo arXiv:1008.3281
\endref 

 \ref\no[BL07]\by  J.~Behrndt and M.~Langer \paper Boundary value
problems for elliptic partial  differential operators on bounded
domains\jour  { J.~Funct.~Anal.}  \vol  243\pages 536--565 \yr2007 \endref

 \ref\no[BL10]\by  J.~Behrndt and M.~Langer \paper Elliptic operators,
 Dirichlet-to-Neumann maps and quasi boundary triples\finalinfo
 preprint, for proceedings of Leiden workshop 2009
\endref

\ref \no[BLLLP10]\by J. Behrndt, M. Langer, I. Lobanov,
V. Lotoreichik and I. Popov  \paper A remark on Schatten-von Neumann
properties of resolvent differences of generalized Robin Laplacians on
bounded domains 
\jour J. Math. Analysis Appl.\vol 371 \yr 2010 \pages 750--758
\endref

\ref\no[B62]\by M. S. Birman\paper Perturbations of the continuous
spectrum of a singular elliptic operator by varying the boundary and
the boundary conditions
\jour Vestnik Leningrad. Univ. \vol 17 \yr 1962 \pages 22--55
 \transl \nofrills English translation in\book
Spectral theory of differential operators, 
 Amer. Math. Soc. Transl. Ser. 2, 225\publ Amer. Math. Soc.\publaddr
Providence, RI \yr 2008 \pages 19--53  
\endref 

\ref\no[BS77] \by   M. S. Birman and M. Z. Solomiak \paper Asymptotics
of the spectrum of pseudo-differential operators with
anisotropic-homogeneous symbols \jour Vestnik Leningrad Univ. \vol 7
\yr 1977 \pages 13--21 \transl\nofrills English translation in \jour Vestnik Leningrad
Univ. Math. \vol 10\yr 1982 \pages 237--247
\endref

\ref\no[BS80] \by   M. S. Birman and M. Z. Solomiak
\paper Asymptotics of the spectrum of variational problems on
solutions of elliptic equations in unbounded domains\jour
  Funkts. Analiz Prilozhen.
  \vol14  \yr1980\pages  27--35\transl\nofrills English translation in
\jour Funct. Anal. Appl. \vol14 \yr1981  \pages267--274
\endref
 

\ref\no[B71]\by 
  L.~Boutet de Monvel  \paper Boundary problems for pseudodifferential
operators\jour  
 {Acta Math.} \vol126\pages  11--51 \yr 1971\endref

\ref\no[BGW09]\by B. M. Brown, G. Grubb, and I. G. Wood \paper $M$-functions for closed
extensions of adjoint pairs of operators with applications to elliptic boundary
problems \jour Math. Nachr. \vol 282\pages 314--347 \yr2009
\endref  

\ref \no[E81] \by G. I. Eskin \book Boundary value problems for
elliptic pseudodifferential equations. Translated from the Russian by
S. Smith. Translations of Mathematical Monographs, 52 \publ American
Mathematical Society \publaddr Providence, R.I \yr 1981 
\endref

\ref\no[F51]
\by Ky Fan
\paper Maximum properties and inequalities for the eigenvalues of
completely continuous operators
\jour Proc. Nat. Acad. Sci. USA
\vol 37
\yr 1951
\pages 760--766
\endref

\ref\no[GM09]\by F. Gesztesy and M. Mitrea \paper Robin-to-Robin maps and
Krein-type resolvent formulas for Schr\"odinger operators on bounded
Lipschitz domains \inbook  Modern Analysis and Applications. The Mark
Krein Centenary Conference, Vol.\ 2. Operator Theory: Advances and Applications  \eds  V.\ Adamyan, Y.\ M.\ Berezansky, I.\
Gohberg, 
M.\ L.\ Gorbachuk, V.\ Gorbachuk, A.\ N.\ Kochubei, H.\ Langer, and G.\
Popov 
\vol 191 \publ Birkh\"auser
\publaddr Basel \yr  2009 \pages 
81--113 
\endref

\ref\no[GM09a]
\by F.~Gesztesy and M.~Mitrea \paper Nonlocal Robin Laplacians and
some remarks on a paper by Filonov on eigenvalue inequalities
\jour J. Diff. Equ. \vol 247 \yr 2009 \pages 2871--2896
\endref

\ref\no[GM10]
\by F.~Gesztesy and M.~Mitrea
\paper A description of all selfadjoint extensions of the Laplacian
 and Krein-type
resolvent formulas in nonsmooth domains\jour J. Analyse Math.\finalinfo
 arXiv:0907.1750, to appear\endref

\ref\no[GK69] \by I. C.  Gohberg and M. G. Krein\book Introduction to the
theory of linear nonselfadjoint operators. Translated from the Russian
by A. Feinstein. Translations of Mathematical Monographs, Vol. 18
\publ American Mathematical Society \publaddr Providence, R.I. \yr
1969 \pages 378  \endref

\ref\no[G68]\by G. Grubb
\paper A characterization of the non-local boundary value problems
associated with an elliptic operator
\jour Ann\. Scuola Norm\. Sup\. Pisa
\vol22
\yr1968\pages425--513
\endref

\ref\no[G70]\by 
{G.~Grubb} \paper Les probl\`emes aux limites g\'en\'eraux d'un
op\'erateur elliptique, provenant de la th\'eorie variationnelle
 \jour{Bull.~ Sc.~Math.} \vol94\pages 113--157 \yr 1970\endref

\ref\no[G74]\by G. Grubb\paper Properties of normal boundary problems for elliptic
even-order systems\jour Ann\. Scuola Norm\. Sup\. Pisa\vol1{\rm
(ser.IV)}\yr1974\pages1--61
\endref

\ref 
\key[G84]
\by G. Grubb
\paper Singular Green operators and their spectral asymptotics
\jour Duke Math. J.
\vol 51
\yr 1984
\pages 477--528
\endref


 \ref\no[G96]\by 
{G.~Grubb}\book Functional Calculus of Pseudodifferential
     Boundary Problems
 Pro\-gress in Math.\ vol.\ 65, Second Edition \publ  Birkh\"auser
\publaddr  Boston \yr 1996\endref

\ref\no[G08]\by G. Grubb \paper Krein resolvent formulas for elliptic
boundary problems in nonsmooth domains \jour Rend. Sem. Mat. Univ.
Pol. Torino \vol 66 \yr2008\pages 13--39
\endref


\ref\no[G09]\by G. Grubb\book Distributions and operators. Graduate
Texts in Mathematics, 252 \publ Springer \publaddr New York\yr 2009
 \endref

\ref\key[G10] \by G. Grubb
\paper
Perturbation of essential spectra of exterior elliptic problems
\finalinfo arXiv:0811.1724, published online 
\jour J. Applicable Analysis
\endref 


\ref\key[G10a]\by G. Grubb \paper Extension theory for elliptic partial
differential operators with pseudodifferential methods 
\finalinfo 
arXiv:1008.1081, for proceedings of Leiden workshop 2009
\endref

\ref\no[GS01]\by G. Grubb and E. Schrohe \paper Trace expansions and the
noncommutative residue for manifolds with boundary \jour J. reine
angew. Math. \vol 536 \yr 2001 \pages 167--207 \endref


\ref\no[J96]\by J.~Johnsen\paper
 Pointwise multiplication of {B}esov and {T}riebel-{L}izorkin spaces
\jour { Math. Nachr.}\vol 175 \pages 85--133 \yr 1995
\endref

\ref\no[L77] \by A. Laptev \paper Spectral asymptotics of a
composition of pseudo-differential operators and reflections from the
boundary
\jour Dokl. Akad. Nauk SSSR \vol 236 \yr 1977\pages 800--830
\transl\nofrills English translation in
\jour Soviet Math. Doklady \vol18 \yr1977  \pages 1273--1276
\endref

\ref\no[L81] \by A. Laptev \paper Spectral asymptotics of a
class of Fourier integral operators
\jour Trudy Mosk. Mat. Obsv. \vol 43 \yr 1981\pages 92--115
\transl\nofrills English translation in
\jour Trans. Moscow Math. Soc.  
\yr1983  \pages 101--127
\endref


\ref\no[LM68]\by  J.-L. Lions and E. Magenes \book  Probl\`emes aux
limites non homog\`enes et applications \vol  1 \publ
 \'Editions Dunod \publaddr Paris \yr 1968
\endref

\ref\no[M10] \by M. M. Malamud \paper Spectral theory of elliptic
operators in exterior domains \jour Russian J. Math. Phys.\vol 17 \yr
2010\pages 96--125 \endref
 

\ref \no[S67] \by R. T. Seeley \paper Complex powers of an elliptic
operator \jour AMS Proc. Symp. Pure Math. \vol 10 \yr 1967 \pages
288--307
\endref



\ref\no[S63] \by E. Shamir\paper Mixed boundary value problems for
elliptic equations in the plane. The $L^{p}$ theory \jour  Ann. Scuola
Norm. Sup. Pisa (3) 17 \yr 1963 \pages 117--139\endref
 

\endRefs
\enddocument